\documentclass{elsart}

\journal{ArXiv}

\usepackage{amsmath}
\usepackage{leftidx}

\begin{document}

\begin{frontmatter}

 \title{Some possible $q$-generalizations of harmonic numbers}
 \author{Istv\'an Mez\H{o}\thanksref{a}}
 \thanks[a]{Present address: University of Debrecen, H-4010, Debrecen, P.O. Box 12, Hungary}
 \address{Department of Applied Mathematics and Probability Theory, Faculty of Informatics, University of Debrecen, Hungary}
 \ead{mezo.istvan@inf.unideb.hu}
 \ead[url]{http://www.inf.unideb.hu/valseg/dolgozok/mezoistvan/mezoistvan.html}

\begin{abstract}We study three different $q$-analogues of the harmonic numbers. As applications, we present some generating functions involving number theoretical functions and give the $q$-generalization of Gosper's exponential generating function of harmonic numbers. We involve also the $q$-gamma and $q$-digamma function.
\end{abstract}

\end{frontmatter}

\newcommand {\Li}{\mathop{\textup{Li}}\nolimits}
\newcommand {\egf}{\mathop{\textup{Egf}}\nolimits}
\newcommand {\gf}{\mathop{\textup{Gf}}\nolimits}
\newtheorem{Proposition}{Proposition}

\section{Introduction}

The harmonic numbers are defined as
\begin{equation}
H_n=\sum_{k=1}^n\frac1k,\quad H_0:=0.\label{hndef}
\end{equation}
Our aim is to find some $q$-analogues of these numbers. We start from some basic identities satisfied by the harmonic numbers. For example,
\begin{eqnarray}
H_n&=&\frac{1}{n!}s(n+1,2),\label{hnst}\\
H_n&=&\sum_{k=1}^n\binom{n}{k}\frac{(-1)^{k+1}}{k},\label{hnbinom}
\end{eqnarray}
where $s(n+1,2)$ is a Stirling number of the first kind \cite{GKP}.

Next we introduce the most basic notions of $q$-calculus. Let
\[[n]_q=\frac{1-q^n}{1-q},\]
and
\[[n]_q!=[n]_q[n-1]_q\cdots[1]_q.\]
Define also $(x;q)_n=(1-x)(1-qx)\cdots(1-q^{n-1}x)$ and $(x;q)_\infty=\lim_{n\rightarrow\infty}(x;q)_n$.
Then the $q$-binomial coefficient with parameter $n$ and $k$ is
\[\binom{n}{k}_q=\frac{[n]_q!}{[k]_q![n-k]_q!}=\frac{(q;q)_n}{(q;q)_k(q;q)_{n-k}}\quad(n\ge k\ge0).\]
According to \eqref{hndef}-\eqref{hnbinom} we define three class of $q$-harmonic numbers as 
\begin{eqnarray}
H_{n,q}^1&=&\sum_{k=1}^n\frac{1}{[k]_q},\\
H_{n,q}^2&=&\frac{1}{[n]_q!}s_q(n+1,2),\label{hndefst}\\
H_{n,q}^3&=&\sum_{k=1}^n\binom{n}{k}_qq^{\binom{k}{2}}\frac{(-1)^{k+1}}{[k]_q},\label{hndefbin1}\\
H_{n,q}^4&=&\ln(q)\sum_{k=1}^n\frac{q^k}{q^k-1}\label{hndeflog}.
\end{eqnarray}
(Let $H_{0,q}^i:=0$ for $i=1,2,3,4$.) Here $s_q(n+1,2)$ is a $q$-Stirling number of the first kind \cite{Charal}.

The first definition appears in \cite{Dilcher,WG}, for example. We point out that the first and second definitions are the same but the third and the fourth are different from each other -- although all of them tend to $H_n$ as $q\rightarrow1$. That $H_{n,q}^4$ differs from the others is obvious because of the presence of $\ln(q)$.

In what follows we deduce some identities and give applications involving these numbers.

\section{Identities involving $H_{n,q}^1$ and $H_{n,q}^2$}

\subsection{Number theoretical results}

In this section we point out that the $H_{n,q}^1$ $q$-harmonic numbers are connected to the divisor function. Moreover, the products of the Riemann zeta function and polylogarithms are Dirichlet generating functions of some interesting number theoretical ``polynomials''.

By definitions,
\[H_{n,q}^1=(1-q)\sum_{k=1}^n\frac{1}{1-q^k}.\]
Cauchy's product gives that
\begin{equation}
\sum_{n\ge 1}H_{n,q}^1x^n=\frac{1-q}{1-x}\sum_{n\ge 1}\frac{x^n}{1-q^n}.\label{gfhnq1}
\end{equation}
We need the notion of Lambert series \cite[p. 257]{HW}. In general, a Lambert series has the form
\[F(q)=\sum_{n\ge 1}a_n\frac{q^n}{1-q^n},\]
where $a_n$ is any suitable sequence.
This connection implies some interesting results. We cite a useful theorem of \cite[Theorem 307]{HW}: if
\[f(s)=\sum_{n\ge 1}\frac{a_n}{n^s},\quad\mbox{and}\quad g(s)=\sum_{n\ge 1}\frac{b_n}{n^s},\]
then
\[F(q)=\sum_{n\ge 1}a_n\frac{q^n}{1-q^n}=\sum_{n\ge 1}b_nq^n\]
holds if and only if
\[\zeta(s)f(s)=g(s),\]
where
\[\zeta(s)=\sum_{n\ge1}\frac{1}{n^s}\]
is the Riemann zeta function.

Let us apply this and \eqref{gfhnq1}:
\begin{equation}
\sum_{n\ge 1}H_{n,q}(qx)^n=\frac{1-q}{1-qx}\sum_{n\ge 1}x^n\frac{q^n}{1-q^n}=\frac{1-q}{1-qx}F(q)\label{hnqs1}
\end{equation}
with $a_n=x^n$. A simple transformation shows that
\begin{equation}
F(q)=\sum_{n\ge 1}q^n\left(\sum_{d|n}x^d\right).\label{FqA}
\end{equation}
Let us introduce the notation
\[A(x,n)=\sum_{d|n}x^d.\]
In special,
\[A(1,n)=d(n)=\sigma_0(n),\]
where $d(n)$ is the divisor function and
\[\sigma_x(n)=\sum_{d|n}d^x.\]
It is worth to realize the symmetry between the power and the base in $A(x,n)$ and $\sigma_x(n)$.

Substituting $x=1$ in \eqref{hnqs1} we have
\begin{Proposition}
\[\sum_{n\ge 1}H_{n,q}q^n=\sum_{n\ge 1}q^nd(n).\]
\end{Proposition}
So the generating function of the divisor function involves the $q$-harmonic numbers. More generally, with the help of \eqref{hnqs1} and \eqref{FqA} and the definition of $A(x,n)$ one can write
\[\sum_{n\ge 1}H_{n,q}(qx)^n=\frac{1-q}{1-qx}\sum_{n\ge 1}q^nA(x,n).\]
Again, by \eqref{hnqs1}
\[\sum_{n\ge 1}x^n\frac{q^n}{1-q^n}=\sum_{n\ge 1}q^nA(x,n).\]
Hence we choose $a_n=x^n$ and $b_n=A(x,n)$ and apply the theorem cited above:
\[\Li_s(x)=\sum_{n\ge 1}\frac{x^n}{n^s},\quad\mbox{and}\quad g(s)=\sum_{n\ge 1}\frac{A(x,n)}{n^s},\]
then
\begin{Proposition}
\[\zeta(s)\Li_s(x)=\sum_{n\ge 1}\frac{A(x,n)}{n^s}.\]
\end{Proposition}
Here $\Li_s(x)$ is the well-known polylogarithm function. In special $\Li_s(1)=\zeta(s)$, so we get the known sum
\[\zeta^2(s)=\sum_{n\ge 1}\frac{d(n)}{n^s}.\]

\subsection{A recursion for $H_{n,q}^1$}

Since
\[H_n=\sum_{k=1}^n\binom{n}{k}\frac{(-1)^{k+1}}{k},\]
we may think that
\begin{equation}
H_{n,q}^1=\sum_{k=1}^n\binom{n}{k}_qa_k\label{hnbak}
\end{equation}
holds for some sequence $a_k$. This is really true but, sadly, $a_k$ does not have a simple form. See the table below.

We shall need the notion of the $q$-Seidel matrix \cite{CHZ}. Given a sequence $a_n$, the $q$-Seidel matrix is associated to the double sequence $a_n^k$ given by the recurrence
\begin{eqnarray*}
a_n^0&=&a_n\quad(n\ge 0),\\
a_n^k&=&q^na_n^{k-1}+a_{n+1}^{k-1}\quad(n\ge 0,k\ge1).
\end{eqnarray*}
In addition, $a_n^0$ is called the initial sequence and $a_0^n$ the final sequence of the $q$-Seidel matrix. Then the identity
\begin{equation}
a_0^n=\sum_{k=0}^n\binom{n}{k}_qa_k^0.\label{Dumont}
\end{equation}
connects the initial and the final sequence.

Define the generating functions of $a_n^0$ and $a_0^n$:
\[a(x)=\sum_{n\ge 0}a_n^0x^n,\quad\overline{a}(x)=\sum_{n\ge 0}a_0^nx^n,\]
and
\[A(x)=\sum_{n\ge 0}a_n^0\frac{x^n}{[n]_q!},\quad\overline{A}(x)=\sum_{n\ge 0}a_0^n\frac{x^n}{[n]_q!}.\]
A proposition given in \cite{CHZ} states that these functions are related by the following equations:
\begin{eqnarray}
\overline{a}(x)&=&\sum_{n\ge 0}a_n^0\frac{x^n}{(x;q)_{n+1}},\label{trf1}\\
\overline{A}(x)&=&e_q(x)A(x),\label{trf2}
\end{eqnarray}
where
\[e_q(x)=\sum_{n\ge 0}\frac{x^n}{[n]_q!}\]
is the $q$-analogue of the exponential function \cite{Gasper}. We introduce the notation $\egf(a_n)$ and $\gf(a_n)$ for the exponential and ordinary generating function of $a_n$, respectively.

To reach our aim posed in \eqref{hnbak}, our approach is as follows. Let the final sequence $b_n=H_{n,q}^1$. We determine the initial sequence $a_n^0=a_n$. Then $\egf(b_n)\equiv\egf(H_{n,q}^1)=e_q\egf(a_n)$. And, to get $\egf(a_n)$ we determine $a_n$ by using \eqref{gfhnq1} and \eqref{trf1}:
\begin{equation}
\gf(b_n)\equiv\gf(H_{n,q}^1)=\frac{1-q}{1-x}\sum_{n\ge 1}\frac{x^n}{1-q^n}=\sum_{n\ge 1}a_n\frac{x^n}{(x;q)_{n+1}}.\label{gfs}
\end{equation}
From this equation $a_n$ can be determined. (Note that $a_0=0$.)

\begin{Proposition}We have
\[H_{n,q}^1=\sum_{k=1}^n\binom{n}{k}_qa_k,\]
where the sequence $a_k$ is determined recursively by
\[\sum_{k=1}^na_kq^{n-k}\binom{n-1}{k-1}_q=\frac{1}{[n]_q}=\frac{1-q}{1-q^n}\quad(a_0:=0).\]

\end{Proposition}

\textit{Proof.} The denominator of the right hand side of \eqref{gfs} is
\begin{equation}
\frac{1}{(x;q)_{n+1}}=\frac{1}{1-x}\frac{1}{(qx;q)_n}=\frac{1}{1-x}\frac{(q^nqx;q)_\infty}{(qx;q)_\infty}.\label{tmp1}
\end{equation}
The $q$-binomial theorem \cite[Section 1.3]{Gasper} states that
\[\frac{(az;q)_\infty}{(z;q)_\infty}=\sum_{k\ge0}\frac{(a;q)_k}{(q;q)_k}z^k.\]
Applying this to \eqref{tmp1},
\[\frac{1}{1-x}\frac{(q^nqx;q)_\infty}{(qx;q)_\infty}=\frac{1}{1-x}\sum_{k\ge0}\frac{(q^n;q)_k}{(q;q)_k}(qx)^k.\]
Thus \eqref{gfs} becomes
\[(1-q)\sum_{n\ge 1}\frac{x^n}{1-q^n}=\sum_{n\ge0}a_nx^n\left(\sum_{k\ge0}\frac{(q^n;q)_k}{(q;q)_k}(qx)^k\right).\]
Let
\[B_{k,n}=\frac{(q^n;q)_k}{(q;q)_k}q^k,\]
for short. Then
\[B_{k,n}=q^k\binom{n+k-1}{k}_q\]
for all $n$ and $k$.
Moreover,
\begin{equation}
(1-q)\sum_{n\ge 1}\frac{x^n}{1-q^n}=\sum_{n\ge0}a_nx^n\left(\sum_{k\ge0}B_{k,n}x^k\right).\label{tmp2}
\end{equation}
If we write the sums term by term, we get
\[a_0(B_{0,0}+B_{1,0}x+B_{2,0}x^2+\cdots)+a_1x^1(B_{0,1}+B_{1,1}x+B_{2,1}x^2+\cdots)+\cdots=\]
\[a_0B_{0,0}+x(a_0B_{1,0}+a_1B_{0,1})+x^2(a_0B_{2,0}+a_1B_{1,1}+a_2B_{0,2})+\cdots\]
Comparing the coefficients here and the left hand side of \eqref{tmp2}, we have
\[\sum_{k=0}^na_kB_{n-k,k}=\frac{1-q}{1-q^n}.\]
Note that -- bacause of \eqref{gfs} -- $a_0$ must be zero. Remember also that $a_k$ is the initial sequence of our $q$-Seidel matrix, so \eqref{Dumont} gives
\begin{equation}
H_{n,q}^1=\sum_{k=1}^n\binom{n}{k}_qa_k.\label{Hnq1an}
\end{equation}
This is our proposition.

\textbf{Remark.} It is worth to present the first terms of the sequence $a_n$:
\begin{eqnarray*}
a_0&=&0,\\
a_1&=&1,\\
a_2&=&-\frac{q^2+q-1}{q+1},\\
a_3&=&\frac{q^5+q^4-q^2-q+1}{q^2+q+1},\\
a_4&=&-\frac{q^9+q^8-2q^5+q^2+q-1}{q^3+q^2+q+1},\\
a_5&=&\frac{q^{14}+q^{13}-q^{10}-q^9-q^8+q^7+q^6+q^5-q^2-q+1}{q^4+q^3+q^2+q+1},
\end{eqnarray*}
\[a_6=-\frac{q^{20}+q^{19}-q^{16}-2q^{14}+q^{12}+q^{11}+q^{10}+q^9-2q^7-q^5+q^2+q-1}{q^5+q^4+q^3+q^2+q+1}.\]
It seemed to be interesting to give a simple formula for the nominator. However, one can easily see that
\[a_k\rightarrow\frac{(-1)^{k+1}}{k}\quad\mbox{as}\quad q\rightarrow 1\quad(k=1,2,3,4,5,6).\]
According to \eqref{hnbinom}, this is plausible for all $a_k$.

As a consequence of \eqref{trf2} and \eqref{Hnq1an}, we have the next connection:
\[\egf(H_{n,q}^1)=e_q\egf(a_n).\]
This is a curious version of Gosper's identity \eqref{Gosper} involving the exponential generating function of the harmonics.

\subsection{The case $H_{n,q}^2$}

To close the case of $H_{n,q}^1$, we remark also that
\begin{equation}
H_{n,q}^1=H_{n,q}^2\label{hn12}
\end{equation}
(see \eqref{hndefst} for the definitions). The $q$-Stirling numbers of the first kind are defined recursively by \cite[p. 155]{Charal}
\begin{equation}
s_q(n+1,k)=s_q(n,k-1)+[n]_qs_q(n,k),\label{sqrec}
\end{equation}
and $s_q(0,0)=1,\;s_q(n,0)=0$ when $n>0$.

These relations imply
\[H_{n,q}^1=\frac{1}{[n]_q!}s_q(n+1,2)=\frac{1}{[n]_q!}s_q(n,1)+\frac{1}{[n-1]_q!}s_q(n,2),\]
hence
\[H_{n,q}^1=H_{n-1,q}^1+\frac{1}{[n]_q!}s_q(n,1)=H_{n-1,q}^1+\frac{[n-1]_q!}{[n]_q!}.\]
Then \eqref{sqrec} gives \eqref{hn12}.

\section{Identities involving $H_{n,q}^3$}

\subsection{A $q$-analogue of Gosper's result}

The exponential generating function of the harmonic numbers is deduced by Gosper \cite{Cvijovic,Mezo}:
\begin{equation}
\sum_{n\ge 1} H_n\frac{x^n}{n!}=xe^x\,\leftidx{_2}{F}{_2}\left(\left.\begin{tabular}{ll}$1$&$1$\\$2$&$2$\end{tabular}\right|-x\right).\label{Gosper}
\end{equation}
Here
\[\leftidx{_r}{F}{_s}\left(\left.\begin{tabular}{llll}$a_1,$&$a_2,$&$\dots,$&$a_r$\\$b_1,$&$b_2$,&$\dots,$&$b_s$\end{tabular}\right|x\right)=\sum_{n\ge 0}\frac{(a_1)_n(a_2)_n\cdots(a_r)_n}{(b_1)_n(b_2)_n\cdots(b_s)_n}\frac{x^n}{n!}\]
is the hypergeometric function with parameters $a_i$ and $b_j$, and
\[(a)_k=a(a+1)(a+2)\cdots(a+k-1)\]
under the agreement $(a)_0=1$.

The $q$-version of the hypergeometric function is called basic hypergeometric function and defined as \cite{Gasper}
\[\leftidx{_r}{\phi}{_s}\left(\left.\begin{tabular}{llll}$a_1,$&$a_2,$&$\dots,$&$a_r$\\$b_1,$&$b_2$,&$\dots,$&$b_s$\end{tabular}\right|q;x\right)=\]
\begin{equation}
=\sum_{n\ge 0}\frac{(a_1;q)_n(a_2;q)_n\cdots(a_r;q)_n}{(b_1;q)_n(b_2;q)_n\cdots(b_s;q)_n}\left((-1)^nq^{\binom{n}{2}}\right)^{1+s-r}\frac{x^n}{(q;q)_n}.\label{bhyper}
\end{equation}
One may see from \eqref{bhyper} that
\[\frac{\mbox{($k+1$)th term}}{\mbox{$k$th term}}=\frac{x(-q^k)^{1+s-r}}{1-q^{k+1}}\frac{(1-a_1q^k)\cdots(1-a_rq^k)}{(1-b_1q^k)\cdots(1-b_sq^k)}.\]

Now we derive the $q$-analogue of Gosper's result \eqref{Gosper} with respecto to $H_{n,q}^3$. Equations \eqref{hndefbin1} and \eqref{trf2} give that
\[\egf(H_{n,q}^3)=e_q\egf\left(q^{\binom{k}{2}}\frac{(-1)^{k+1}}{[k]_q}\right).\]
Hence it is enough to determine the sum
\[\sum_{k\ge 1}q^{\binom{k}{2}}\frac{(-1)^{k+1}}{[k]_q}\frac{x^k}{[k]_q!}.\]
We would like to express this function as a basic hypergeometric series. Thus $k$ should run from $0$, so we write
\[x\sum_{k\ge 0}q^{\binom{k+1}{2}}\frac{(-1)^{k+2}}{[k+1]_q}\frac{x^k}{[k+1]_q!},\]
whence
\[\frac{\mbox{($k+1$)th term}}{\mbox{$k$th term}}=(1-q)\frac{(-xq)q^k}{1-q^{k+1}}\frac{(1-qq^k)^2}{(1-q^2q^k)^2}.\]
Since $s=r=2$ in our case, the place of $(-1)^n$ is indifferent. Finally we can state the following
\begin{Proposition}We have
\[\sum_{n\ge 1}H_{n,q}^3\frac{x^n}{[n]_q!}=(1-q)xe_q(x)\leftidx{_2}{\phi}{_2}\left(\left.\begin{tabular}{llll}$q,$&$q$\\$q^2,$&$q^2$\end{tabular}\right|q;-qx\right).\]
\end{Proposition}

\subsection{The Hockey Stick Theorem -- a $q$-analogue}

It is natural to ask, what is the recursion satisfied by $H_{n,q}^3$? If $H_{n,q}^3=H_{n-1,q}^3+\frac{1}{[n]_q}$ would be true, we knew that $H_{n,q}^3=H_{n,q}^1$. But this relation does not hold. In fact, the next recursion is valid.

\begin{Proposition}\label{PropHnq3}For $H_{n,q}^3$
\[H_{n,q}^3=H_{n-1,q}^3+q-\frac{[n-1]_q}{[n]_q}.\]
\end{Proposition}
We see that the limit of $q-\frac{[n-1]_q}{[n]_q}$ is $1-\frac{n-1}{n}=\frac1n$. So the standard recursion for harmonic numbers holds just asymptotically. To prove the proposition, we need the following statement. The standard (not $q$-) version can be found at the webpage http://binomial.csueastbay.edu/IdentitiesNamed.htm, Catalog \#: 3400001.

\begin{Proposition}[Hockey Stick Theorem]
\[\sum_{k=1}^{n-j}(-1)^{k+1}\binom{n}{j+k}=\binom{n-1}{j}.\]
\end{Proposition}

The name comes from the fact that the summands and the sum has a shape in the Pascal triangle like a hockey stick.

We prove the $q$-analogue version.
\begin{Proposition}[Hockey Stick Theorem -- a $q$-analogue]
\[\sum_{k=1}^{n-j}(-1)^{k+1}q^{\binom{j+k}{2}}\binom{n}{j+k}_q=\binom{n-1}{j}_q.\]
\end{Proposition}

\textit{Proof.} Write the sum term by term:
\[q^{\binom{j+1}{2}}\binom{n}{j+1}_q-q^{\binom{j+2}{2}}\binom{n}{j+2}_q+q^{\binom{j+3}{2}}\binom{n}{j+3}_q-\cdots+\]
\[(-1)^{n-j}q^{\binom{n-1}{2}}\binom{n}{n-1}_q+(-1)^{n-j+1}q^{\binom{n}{2}}\binom{n}{n}_q.\]
The binomial coefficients are rewritten with the recursion \cite{Gasper}
\[\binom{n}{k}_q=q^k\binom{n-1}{k}_q+\binom{n-1}{k-1}_q:\]
\[q^{\binom{j+1}{2}}\left[q^{j+1}\binom{n-1}{j+1}_q+\binom{n-1}{j}_q\right]-q^{\binom{j+2}{2}}\left[q^{j+2}\binom{n-1}{j+2}_q+\binom{n-1}{j+1}_q\right]+\]
\[q^{\binom{j+3}{2}}\left[q^{j+3}\binom{n-1}{j+3}_q+\binom{n-1}{j+2}_q\right]-\cdots+\]
\[(-1)^{n-j}q^{\binom{n-1}{2}}\left[q^{n-1}\binom{n-1}{n-1}_q+\binom{n-1}{n-2}_q\right]+(-1)^{n-j+1}q^{\binom{n}{2}}\left[0+\binom{n-1}{n-1}_q\right].\]
Realize that the first members in the square brackets with the $q$-power coefficients are cancelled by the second member in the next square bracket:
\[q^{\binom{j+k}{2}}q^{j+k}\binom{n-1}{j+k}_q-q^{\binom{j+k+1}{2}}\binom{n-1}{j+k}_q=0.\]
This is true for all $k=1,2,\dots,n-j-1$.

Now we are ready to prove Proposition \ref{PropHnq3}. A consequence of \eqref{hndefbin1}:
\[H_{n,q}^3-H_{n-1,q}^3=\left(\binom{n}{1}_q-\binom{n-1}{1}_q\right)\frac{q^{\binom{1}{2}}}{[1]_q}+\left(-\binom{n}{2}_q+\binom{n-1}{2}_q\right)\frac{q^{\binom{2}{2}}}{[2]_q}+\]
\[\left(\binom{n}{3}_q-\binom{n-1}{3}_q\right)\frac{q^{\binom{3}{2}}}{[3]_q}+\cdots=\]
\[\left(\binom{n-1}{0}_q-(1-q)\binom{n-1}{1}_q\right)\frac{q^{\binom{1}{2}}}{[1]_q}+\left(-\binom{n-1}{1}_q+(1-q^2)\binom{n-1}{2}_q\right)\frac{q^{\binom{2}{2}}}{[2]_q}+\]
\[\left(\binom{n-1}{2}_q-(1-q^3)\binom{n-1}{3}_q\right)\frac{q^{\binom{3}{2}}}{[3]_q}+\cdots=\]
\[1-(1-q)\binom{n-1}{1}_q+\sum_{k=1}^{n-1}\binom{n-1}{k}_q\frac{(-1)^k}{[k+1]_q}q^{\binom{k+1}{2}}+\sum_{k=2}^{n-1}\binom{n-1}{k}_q\frac{(-1)^k(1-q^k)}{[k]_q}q^{\binom{k}{2}}=\]
\[1+\sum_{k=1}^{n-1}\binom{n-1}{k}_q\frac{(-1)^k}{[k+1]_q}q^{\binom{k+1}{2}}+\sum_{k=1}^{n-1}\binom{n-1}{k}_q\frac{(-1)^k(1-q^k)}{[k]_q}q^{\binom{k}{2}}=\]
\[1+\frac{1}{[n]_q}\sum_{k=1}^{n-1}\binom{n}{k+1}_q(-1)^kq^{\binom{k+1}{2}}+(1-q)\sum_{k=1}^{n-1}\binom{n-1}{k}_q(-1)^kq^{\binom{k}{2}}.\]
The first sum is exactly the left hand side of the $q$-Hockey Stick Theorem with $j=1$, up to a minus sign. The second sum equals $-1$ \cite[p. 153]{Charal}. Thus
\[H_{n,q}^3-H_{n-1,q}^3=1-\frac{1}{[n]_q}\binom{n-1}{1}_q+q-1.\]
So we have the desired result.

\section{Identity involving $H_{n,q}^4$}

The digamma function is defined as
\[\psi(x)=\frac{\Gamma'(x)}{\Gamma(x)},\]
where $\Gamma(x)$ is the Euler gamma function. Latter satisfies the functional equation
\begin{equation}
\Gamma(x+1)=x\Gamma(x).\label{gammafe}
\end{equation}
The harmonic numbers are connected to the digamma function (as one immediately sees by the logarithmic derivative of \eqref{gammafe}):
\[H_n=\psi(n+1)+\gamma.\]
Here $\gamma=-\psi(1)$ is the Euler-Mascheroni constant. Our goal is to find the $q$-analogue of this formula.

Let us start from the definition of the $q$-gamma function (see \cite{Gasper}):
\[\Gamma_q(x)=\frac{(q;q)_\infty}{(q^x;q)_\infty}(1-q)^{1-x}.\]
The $q$-digamma function is simply the logarithmic derivative of $\Gamma_q(x)$ as in the ordinary case (\cite{Gasper,Jackson,Thomae}):
\[\psi_q(x)=-\ln(1-q)+\ln(q)\sum_{n\ge 0}\frac{q^{n+x}}{1-q^{n+x}}.\]

It is known that $\Gamma_q(x)$ satisfies the next $q$-version of $\eqref{gammafe}$:
\[\Gamma_q(x+1)=[x]_q\Gamma_q(x).\]
Then taking logarithm and then the derivative, we get
\[\psi_q(x+1)=\frac{d}{dx}\ln([x]_q)+\psi_q(x).\]
Since
\[\frac{d}{dx}\ln([x]_q)=\frac{\ln(q)q^x}{q^x-1},\]
the previous equation is rewritten recursively as
\[\psi_q(x+1)=\frac{\ln(q)q^x}{q^x-1}+\frac{\ln(q)q^{x-1}}{q^{x-1}-1}+\cdots+\frac{\ln(q)q}{q-1}+\psi_q(1).\]
So definition \eqref{hndeflog} seems to be correct.

On the other hand, we may define the $q$-analogue of the Euler-Mascheroni constant as
\[\gamma_q=-\psi_q(1)=\ln(1-q)-\ln(q)\sum_{n\ge 1}\frac{q^n}{1-q^n},\]
or more simply,
\[\gamma_q=\ln(1-q)+\frac{\ln(q)}{q-1}\sum_{n\ge 1}\frac{q^n}{[n]_q}.\]

With these we have
\[H_{n,q}^4=\psi_q(n+1)+\gamma_q.\]


\begin{thebibliography}{AAAAAA}
\bibitem[1]{Charal}Ch. Charalambides, Combinatorial Methods in Discrete Dis\-tri\-bu\-ti\-ons, Wi\-ley\,--\,In\-ter\-sci\-en\-ce, 2005.
\bibitem[2]{CHZ}R. J. Clarke, G.-N. Han and J. Zeng, A combinatorial interpretation of the Seidel generation of $q$-derangement numbers, Annals of Comb. 4 (1997), 313-327.
\bibitem[3]{Cvijovic}D. Cvijovi\'c, The Dattoli-Srivastava conjectures concerning generating functions involving the harmonic numbers, Appl. Math. Comput. 215(11) (2010), 4040-4043.
\bibitem[4]{Dilcher}K. Dilcher, Determinant expressions for $q$-harmonic congruences and degenerate Bernoulli numbers, The Electronic Journal of Combinatorics 15 (2008), \#R63.
\bibitem[5]{Gasper}G. Gasper and M. Rahman, Basic Hypergeometric Series (second edition), Cambridge University Press, 2004.
\bibitem[6]{GKP}R. L. Graham, D. E. Knuth and O. Patashnik, \emph{Concrete Mathematics}, Addison Wesley, 1993.
\bibitem[7]{HW}G. H. Hardy, E. M. Wright, An Introduction to the Theory of Numbers, Oxford University Press, 1979.
\bibitem[8]{Jackson}F. H. Jackson, A generalization of the functions $\Gamma(n)$ and $x^n$, Proc. Roy. Soc. London, 74 (1904), 64-72.
\bibitem[9]{Mezo}I. Mez\H{o} and A. Dil, Euler-Seidel method for certain combinatorial numbers and a new characterization of Fibonacci sequence, Cent. Eur. J. Math. 7(2) (2009), 310-321.
\bibitem[10]{Thomae}J. Thomae, Beitr\"age zur Theorie der durch die Heinesche Reihe\dots, J. Reine und Angew. Math. 70 (1869), 258-281.
\bibitem[11]{WG}C. Wei, Q. Gu, $q$-generalizations of a family of harmonic number identities, Adv. Appl. Math. 45(1) (2010), 24-27.
\end{thebibliography}
\end{document}